\documentclass{amsart}
\usepackage{latexsym}
\usepackage{amsfonts}
\usepackage{amssymb}

\newtheorem{theorem}{Theorem}[section]
\newtheorem{corollary}[theorem]{Corollary}
\newtheorem{lemma}[theorem]{Lemma}
\newtheorem{proposition}[theorem]{Proposition}

\theoremstyle{remark}
\newtheorem{remark}{Remark}[section]

\numberwithin{equation}{section}

\newcommand{\g}{\mathfrak{g}}
\newcommand{\ak}{\mathfrak{k}}
\newcommand{\at}{\mathfrak{t}}
\newcommand{\am}{\mathfrak{m}}
\newcommand{\ZZ}{\mathbb{Z}}

\newcommand{\Zn}{\mathbb{Z}/n\mathbb{Z}}
\newcommand{\Ad}{\mathrm{Ad}}
\newcommand{\Inn}{\mathrm{Inn}}
\newcommand{\rank}{\mathrm{rank}}
\newcommand{\diag}{\mathrm{diag}}

\begin{document}

\title{Twisted Weyl groups of Lie groups and nonabelian cohomology}

\author{Jinpeng An}
\address{Department of Mathematics, ETH Zurich, 8092 Zurich, Switzerland}
\email{anjp@math.ethz.ch}



\begin{abstract}
For a cyclic group $A$ and a connected Lie group $G$ with an
$A$-module structure (with the additional conditions that $G$ is
compact and the $A$-module structure on $G$ is $1$-semisimple if
$A\cong\ZZ$), we define the twisted Weyl group $W=W(G,A,T)$, which
acts on $T$ and $H^1(A,T)$, where $T$ is a maximal compact torus
of $G_0^A$, the identity component of the group of invariants
$G^A$. We then prove that the natural map $W\backslash
H^1(A,T)\rightarrow H^1(A,G)$ is a bijection, reducing the
calculation of $H^1(A,G)$ to the calculation of the action of $W$
on $T$. We also prove some properties of the twisted Weyl group
$W$, one of which is that $W$ is a finite group. A new proof of a
known result concerning the ranks of groups of invariants with
respect to automorphisms of a compact Lie group is also given.
\end{abstract}

\maketitle

%%%%%%%%%%%%%%%%%%%%%%%%%%%%%%%%%%%%%%%%%%%%%%%%%%%%%%%%%%% Section 1

\section{Introduction}
Let $A$ be a group, and let $G$ be a connected Lie group with an
$A$-module structure. It is proved in \cite{AW} that, if $A$ is
finite, then there exists a maximal compact subgroup $K$ of $G$
which is an $A$-submodule of $G$, and for every such $K$,
$H^1(A,K)\rightarrow H^1(A,G)$ is a bijection. This reduces the
calculation of $H^1(A,G)$ to the calculation of $H^1(A,K)$ for $K$
compact.

Based on some results in \cite{AW}, in this paper we go further
along this direction. We reduce the calculation of $H^1(\ZZ,G)$
for $G$ compact, and $H^1(\Zn,G)$ for $G$ general, to the
calculation of the action of the twisted Weyl group on a compact
torus of $G$.

Let $A$ be a cyclic group, that is, $A$ is isomorphic to $\ZZ$ or
$\Zn$, and let $G$ be a connected Lie group with an $A$-module
structure. In the case of $A\cong\ZZ$, we always assume that $G$
is compact and the $\ZZ$-module structure on $G$ is
$1$-semisimple, that is,
$$\ker(1-d\sigma)=\ker((1-d\sigma)^2)$$ for a generator $\sigma$ of
$A$. Let $T$ be a maximal compact torus of $G_0^A$, the identity
component of $G^A=\{g\in G|a(g)=g, \forall a\in A\}$. Then it is
proved in \cite{AW} that the natural map $H^1(A,T)\rightarrow
H^1(A,G)$ is surjective (see Theorems 4.1 and 5.1 in \cite{AW}).

Let $\sigma$ be a generator of $A$. Denote
$$Z_\sigma(T)=\{g\in G|gt\sigma(g)^{-1}=t,\forall t\in T\},$$
$$N_\sigma(T)=\{g\in G|gT\sigma(g)^{-1}=T\},$$
then $Z_\sigma(T)$ and $N_\sigma(T)$ are closed subgroups of $G$,
$Z_\sigma(T)$ is normal in $N_\sigma(T)$. Define the \emph{twisted
Weyl group} $W=W(G,A,T)$ by
$$W=N_\sigma(T)/Z_\sigma(T).$$
Then it can be proved that $W$ is independent of the choice of
$\sigma$, and, as an abstract group, $W$ is independent of the
choice of $T$. It can also be proved that $W$ is a finite group.

Since $T$ is abelian and $A$ acts trivially on $T$, $H^1(A,T)$ may
be identified with the set of cocycles $Z^1(A,T)$, which consists
of homomorphisms $A\rightarrow T$. The twisted Weyl group $W$ acts
on $Z^1(A,T)\cong H^1(A,T)$ by
$$(w.\alpha)(a)=g\alpha(a)a(g)^{-1},$$ where $w\in W,
\alpha\in Z^1(A,T), a\in A, g\in w$. Denote by $W\backslash
H^1(A,T)$ the space of $W$-orbits in $H^1(A,T)$. Then the natural
map $H^1(A,T)\rightarrow H^1(A,G)$ reduces to a map $W\backslash
H^1(A,T)\rightarrow H^1(A,G)$.

The main result of this paper is the following assertion.

\begin{theorem}\label{T:main}
Suppose $A$ is a cyclic group, $G$ is a connected Lie group with
an $A$-module structure. If $A\cong\ZZ$, we assume that $G$ is
compact and the $\ZZ$-module structure on $G$ is $1$-semisimple.
Let $T$ be a maximal compact torus of $G_0^A$, $W=W(G,A,T)$ the
associated twisted Weyl group. Then the map $W\backslash
H^1(A,T)\rightarrow H^1(A,G)$ is a bijection.
\end{theorem}

Upon choosing a generator of $A$, we may identify $H^1(A,T)$ with
$T$ (if $A\cong\ZZ$) or the finite subgroup $E_n(T)=\{t\in
T|t^n=e\}$ of $T$ (if $A\cong\Zn$). So Theorem \ref{T:main}
reduces the calculation of $H^1(A,G)$ to the calculation of the
action of $W$ on $T$ or $E_n(T)$. (This action will be defined
directly in Section 3.)

In the case of $A\cong\ZZ/2\ZZ$, $G$ is a complex reductive group,
and $\ZZ/2\ZZ$ acts on $G$ by complex conjugation with respect to
a real form of $G$, the bijectivity of the map $W\backslash
H^1(\ZZ/2\ZZ,T)\rightarrow H^1(\ZZ/2\ZZ,G)$ was proved by Borovoi
\cite{Bo}.

Our proof of Theorem \ref{T:main} relies on the notion of twisted
conjugate actions of Lie groups associated with automorphisms.
Recall that for a Lie group $G$ with an automorphism $\sigma$, the
twisted conjugate action of $G$ on itself associated with $\sigma$
is defined by
$$\tau_g(h)=gh\sigma(g)^{-1}.$$ Two elements $g_1,g_2\in G$ are
\emph{$\sigma$-conjugate} if they lie in the same orbit of the
twisted conjugate action associated with $\sigma$. The utility of
this notion is based on the following fact: If $A$ is cyclic with
a generator $\sigma$, then two cocycles $\alpha_1, \alpha_2\in
Z^1(A,G)$ are cohomologous if and only if $\alpha_1(\sigma)$ and
$\alpha_2(\sigma)$ are $\sigma$-conjugate.

Similar to the definition of $W(G,A,T)$, we can define the twisted
Weyl group $W(G,\sigma,T)$ associated with an individual
automorphism $\sigma$ of $G$, which acts naturally on the maximal
compact torus $T$ of $G_0^\sigma$. The proof of Theorem
\ref{T:main} is based on the following result.

\begin{theorem}\label{T:1.2}
Let $G$ be a connected compact Lie group with a $1$-semisimple
automorphism $\sigma$, $T$ a maximal torus of $G_0^\sigma$. Then
two elements $t_1,t_2\in T$ are $\sigma$-conjugate if and only if
they lie in the same $W(G,\sigma,T)$-orbit.
\end{theorem}

To prove Theorem \ref{T:1.2}, we need the following known result.

\begin{theorem}\label{T:1.3}
Let $\sigma_1,\sigma_2$ be automorphisms of a connected compact
Lie group $G$. If $\sigma_2\circ\sigma_1^{-1}$ is an inner
automorphism of $G$, then $\rank G_0^{\sigma_1}=\rank
G_0^{\sigma_2}$.
\end{theorem}

Although the proof of Theorem \ref{T:1.3} is implicitly included
in Gantmacher \cite{Ga}, the author can not find an explicit
reference for it. Using properties of twisted conjugate actions,
we give a new proof of Theorem \ref{T:1.3} in Section 2. In
Section 3 we define the twisted Weyl group of a Lie group
associated with an automorphism, prove some of its properties, and
give the proof of Theorem \ref{T:1.2}. The proof of Theorem
\ref{T:main} is given is Section 4.

The author would like to thank Jiu-Kang Yu for kind guidance.

%%%%%%%%%%%%%%%%%%%%%%%%%%%%%%%%%%%%%%%%%%%%%%%%%%%%%%%%%%% Section 2

\section{A preliminary result}

We discuss in this section a known result concerning groups of
invariants with respect to automorphisms of a compact Lie group,
which will be used in the proofs of some results in later
sections.

For a connected Lie group $G$ with an automorphism $\sigma$, we
denote by $G^\sigma=\{g\in G|\sigma(g)=g\}$ the group of
invariants, and denote by $G_0^\sigma$ the identity component of
$G^\sigma$.

\begin{theorem}\label{T:application}
Let $\sigma_1,\sigma_2$ be automorphisms of a connected compact
Lie group $G$. If $\sigma_2\circ\sigma_1^{-1}$ is an inner
automorphism, then $\rank G_0^{\sigma_1}=\rank G_0^{\sigma_2}$.
\end{theorem}

Theorem \ref{T:application} was implicitly proved in \cite{Ga}
(see also \cite{OV}, Chapter 4, Section 4). But the author can not
find an explicit reference. Here we give a new proof of this
result, with the aid of some properties of twisted conjugate
actions. This may also be viewed as an application of twisted
conjugate actions.

We first prove the following basic result.

\begin{lemma}\label{L:dimension}
Suppose a Lie group $G$ acts smoothly on a smooth manifold $M$.
Let $n=\dim M$, $d$ be the maximal dimension of $G$-orbits in $M$.
If $N$ is a closed submanifold of $M$ such that the intersection
of every $G$-orbit with $N$ is a nonempty and discrete subset of
$M$, then $\dim N=m-d$.
\end{lemma}

\begin{proof}
Denote the action of $G$ on $M$ by $\varphi:G\times M\rightarrow
M$. For $x\in M$, denote the $G$-orbit containing $x$ by $O_x$.
Let $M'=\{x\in M|\dim O_x=d\}$. Then for a point $x\in M$, $x\in
M'$ if and only if the differential of the map
$\varphi(\cdot,x):G\rightarrow M$ has maximal rank. So $M\setminus
M'$ can be locally described as the zero locus of a smooth
function on $M$. Then $M'$ is open in $M$.

Since $\varphi$ is transversal to $M$, $\varphi^{-1}(N)$ is a
closed submanifold of $G\times M$. Denote by $\pi_2:G\times
M\rightarrow M$ the projection to the second factor, and consider
the restriction
$\pi_2|_{\varphi^{-1}(N)}:\varphi^{-1}(N)\rightarrow M$ of $\pi_2$
to $\varphi^{-1}(N)$. By Sard's Theorem, we can choose a regular
value $p\in M'$ of $\pi_2|_{\varphi^{-1}(N)}$. Since $O_p\cap N$
is nonempty, there exists $g\in G$ such that $q=\varphi(g,p)\in
N$. Since $(g,p)$ is a regular point of
$\pi_2|_{\varphi^{-1}(N)}$,
$T_{(g,p)}\varphi^{-1}(N)+T_{(g,p)}G\times\{p\}=T_{(g,p)}G\times
M$. This implies $T_qN+T_qO_p=T_qM$. So $N$ intersects $O_p$
transversally at $q$. But $N\cap O_p$ is discrete, so $\dim N+\dim
O_p=n$. Note that $p\in M'$ implies $\dim O_p=d$, we have $\dim
N=n-d$.
\end{proof}

\begin{proof} [Proof of Theorem \ref{T:application}]
Let $\tau_1, \tau_2$ be the twisted conjugate actions associated
with $\sigma_1$ and $\sigma_2$, respectively. We first assume that
$G$ is semisimple. Let $T_1, T_2$ be maximal tori of
$G_0^{\sigma_1}$ and $G_0^{\sigma_2}$, respectively. Then by
Theorem 6.1 in \cite{AW}, the intersection of every $\tau_1$-orbit
(resp. $\tau_2$-orbit) with $T_1$ (resp. $T_2$) is nonempty and
finite. By Lemma \ref{L:dimension}, we have
$$\dim T_1=\dim G-\text{the maximal dimension of $\tau_1$-orbits},$$
$$\dim T_2=\dim G-\text{the maximal dimension of $\tau_2$-orbits}.$$
But by Proposition 6.1 in \cite{AW}, $\tau_1$ and $\tau_2$ are
equivalent. So the maximal dimension of $\tau_1$-orbits and that
of $\tau_2$-orbits coincide. Thus we have $\dim T_1=\dim T_2$,
that is, $\rank G_0^{\sigma_1}=\rank G_0^{\sigma_2}$.

For the general case, let $G_s$ be the semisimple part of $G$, and
let $G_t$ be the identity component of the center of $G$, which is
a compact torus. Then $G_s$ and $G_t$ are invariant under
$\sigma_1$ and $\sigma_2$. Since $\sigma_2\circ\sigma_1^{-1}$ is
inner, $(G_t)_0^{\sigma_1}=(G_t)_0^{\sigma_2}$. But we have proved
that $\rank (G_s)_0^{\sigma_1}=\rank (G_s)_0^{\sigma_2}$. So
$$\rank G_0^{\sigma_1}=\rank (G_s)_0^{\sigma_1}+\dim
(G_t)_0^{\sigma_1}=\rank (G_s)_0^{\sigma_2}+\dim
(G_t)_0^{\sigma_2}=\rank G_0^{\sigma_2}.$$
\end{proof}

%%%%%%%%%%%%%%%%%%%%%%%%%%%%%%%%%%%%%%%%%%%%%%%%%%%%%%%%%%% Section 3

\section{Twisted Weyl groups associated with automorphisms}

In order to define the twisted Weyl group associated with an
$A$-module structure on a Lie group, in this section we first
consider the twisted Weyl group associated with a single
automorphism of a Lie group.

Let $G$ be a connected Lie group with an automorphism $\sigma$.
For technical reasons, we always make the following assumptions.\\
(1) \emph{If $G$ is compact, we assume that $\sigma$ is
$1$-semisimple, that is,}
$$\ker(1-d\sigma)=\ker((1-d\sigma)^2);$$ (2) \emph{If $G$ is
noncompact, we assume that $\sigma$ is of
finite order.}\\
Note that an automorphism of finite order is $1$-semisimple. Let
$\tau$ be the twisted conjugate action of $G$ on itself associated
with $\sigma$, which is defined by $\tau_g(h)=gh\sigma(g)^{-1}$.
For a closed subgroup $H$ of $G$, denote
$$Z_\sigma(H)=Z_{\sigma,G}(H)=\{g\in G|\tau_g(t)=t,\forall t\in H\},$$
$$N_\sigma(H)=N_{\sigma,G}(H)=\{g\in G|\tau_g(H)=H\}.$$
Then $Z_\sigma(H)$ and $N_\sigma(H)$ are closed subgroups of $G$,
and $Z_\sigma(H)$ is normal in $Z_\sigma(H)$. Now let $T$ be a
maximal compact torus of $G_0^\sigma$. Define the \emph{twisted
Weyl group associated with $\sigma$} by
$$W(G,\sigma,T)=N_\sigma(T)/Z_\sigma(T).$$ We will denote
$W(G,\sigma,T)$ by $W(\sigma,T)$, or $W(T)$, or simply $W$, if the
omitted data are explicit from the context.

\begin{proposition}\label{P:independent}
As an abstract group, $W(T)$ is independent of the choices of $T$.
\end{proposition}

\begin{proof}
Let $T'$ be another maximal compact torus of $G_0^\sigma$. By
Proposition 5.1 in \cite{AW}, there exists $g\in G_0^\sigma$ such
that $T'=gTg^{-1}$. It is easily verified that
$Z_\sigma(T')=gZ_\sigma(T)g^{-1}$,
$N_\sigma(T')=gN_\sigma(T)g^{-1}$. So $W(T)\cong W(T')$.
\end{proof}

\begin{proposition}\label{P:lie algebra}
(i) $Z_\sigma(T)\subset G^\sigma$;\\
(ii) The Lie algebras of $Z_\sigma(T)$ and $N_\sigma(T)$
coincide;\\
(iii) If $G$ is compact, then the Lie algebras of $Z_\sigma(T)$
and $N_\sigma(T)$ coincide with the Lie algebra of $T$.
\end{proposition}

\begin{proof}
(i) Let $g\in Z_\sigma(T)$, then for every $t\in T$,
$\tau_g(t)=t$. In particular, $\tau_g(e)=g\sigma(g)^{-1}=e$, that
is, $\sigma(g)=g$. Hence $g\in G^\sigma$.

(ii) If $X$ belongs to the Lie algebra of $N_\sigma(T)$, then for
every $s\in\mathbb{R}$ and every $t\in T$,
$\tau_{e^{sX}}(t)=e^{sX}te^{-sd\sigma(X)}\in T$. This implies
$(1-\Ad(t)d\sigma)X\in\at$, the Lie algebra of $T$. Since
$\at\subset\ker(1-\Ad(t)d\sigma)$, $(1-\Ad(t)d\sigma)^2X=0$. Since
the mutually commutative endomorphisms $d\sigma$ and $\Ad(t)$ are
$1$-semisimple, $\Ad(t)d\sigma$ is $1$-semisimple. So we have
$(1-\Ad(t)d\sigma)X=0$ for every $t\in T$, that is, $X$ belongs to
the Lie algebra of $Z_\sigma(T)$.

(iii) Suppose $G$ is compact. By (i), the identity component
$Z_\sigma(T)_0$ of $Z_\sigma(T)$ is contained in $G_0^\sigma$. So
$Z_\sigma(T)_0\subset Z_\sigma(T)\cap
G_0^\sigma=Z_{G_0^\sigma}(T)=T$, where $Z_{G_0^\sigma}(T)$ is the
centralizer of $T$ in $G_0^\sigma$. But $T\subset Z_\sigma(T)_0$,
so $Z_\sigma(T)_0=T$. Hence the Lie algebra of $Z_\sigma(T)$
coincides with the Lie algebra of $T$, and then (iii) follows from
(ii).
\end{proof}

\begin{remark}
Suppose $G$ is compact. If $\sigma$ is the identity, we have
$Z_\sigma(T)=Z_G(T)=T$. But this does not hold for general
$\sigma$, although $Z_\sigma(T)$ and $T$ have the same Lie
algebra. For example, let $G=U(3)$ with automorphism $\sigma$
defined by $\sigma(g)=\overline{g}$. Then $G^\sigma=O(3)$,
$G^\sigma_0=SO(3)$. $T=\diag(SO(2),1)$ is a maximal torus of
$G_0^\sigma$, $g=\diag(1,1,-1)\in Z_\sigma(T)$, but $g\notin T$.
\end{remark}

\begin{remark}
It is not necessary that $N_\sigma(T)\subset G^\sigma$. Consider
also the example in the above remark. Then $h=\diag(i,i,1)\in
N_\sigma(T)$, but $h\notin G^\sigma$.
\end{remark}

By (ii) of Proposition \ref{P:lie algebra}, the twisted Weyl group
$W$ is finite if $G$ is compact. We claim that it is also finite
in the noncompact case. To prove this, we need some preliminaries,
some of which are also used in the next section.

Suppose $G$ is a connected Lie group with an automorphism $\sigma$
of finite order. By Theorem 3.1 in \cite{AW}, there always exists
a maximal compact subgroup $K$ of $G$ which is $\sigma$-invariant.
We first prove two lemmas.

\begin{lemma}\label{L:basic}
Let $G$ be a connected Lie group with an automorphism $\sigma$ of
finite order, and let $K$ be a $\sigma$-invariant maximal compact
subgroup of $G$. Then every $g\in G$ admits a decomposition $g=kp$
such that $k\in K$, and such that for every $k_1,k_2\in K$, if
$\tau_g(k_1)=k_2$, then $\tau_k(k_1)=k_2$, $\tau_p(k_1)=k_1$.
\end{lemma}

\begin{proof}
Let $\g$ and $\ak$ be the Lie algebras of $G$ and $K$,
respectively. By \cite{Ho}, Chapter XV, Theorem 3.1, there exist
linear subspaces $\am_1,\cdots,\am_r$ of $\g$ with
$\g=\ak\oplus\am_1\oplus\cdots\oplus\am_r$ such that
$\Ad(k)(\am_i)=\am_i, \forall k\in K, i\in\{1,\cdots,r\}$, and
such that the map
$\varphi:K\times\am_1\times\cdots\times\am_r\rightarrow G$ defined
by $\varphi(k,X_1,\cdots,X_r)=ke^{X_1}\cdots e^{X_r}$ is a
diffeomorphism. For $g\in G$, write $g$ as $g=ke^{X_1}\cdots
e^{X_r}$, where $k\in K$, $X_i\in\am_i$. Let $k_1, k_2\in K$ such
that $\tau_g(k_1)=gk_1\sigma(g)^{-1}=k_2$. Rewrite this equality
as $k_1^{-1}gk_1=k_1^{-1}k_2\sigma(g)$, then we have
$$
(k_1^{-1}kk_1)e^{\Ad(k_1^{-1})X_1}\cdots
e^{\Ad(k_1^{-1})X_r}=(k_1^{-1}k_2\sigma(k))e^{d\sigma(X_1)}\cdots
e^{d\sigma(X_r)},
$$
that is,
$$
\varphi(k_1^{-1}kk_1,\Ad(k_1^{-1})X_1,\cdots,\Ad(k_1^{-1})X_r)=
\varphi(k_1^{-1}k_2\sigma(k),d\sigma(X_1),\cdots,d\sigma(X_r)).
$$
Since $\varphi$ is a diffeomorphism, we have
$k_1^{-1}kk_1=k_1^{-1}k_2\sigma(k)$, that is, $\tau_k(k_1)=k_2$.
Let $p=e^{X_1}\cdots e^{X_r}=k^{-1}g$, then
$\tau_p(k_1)=\tau_{k^{-1}}\circ\tau_g(k_1)=\tau_{k^{-1}}(k_2)=k_1$.
This proves the lemma.
\end{proof}

\begin{lemma}\label{L:normalizer}
Let $G$ be a connected Lie group with an automorphism $\sigma$ of
finite order, and let $K$ be a $\sigma$-invariant maximal compact
torus of $G$. Then for every closed subgroup $H$ of $K$,
$N_{\sigma,G}(H)=N_{\sigma,K}(H)\cdot Z_{\sigma,G}(H)$.
\end{lemma}

\begin{proof}
It is obvious that $N_{\sigma,K}(H)\cdot Z_{\sigma,G}(H)\subset
N_{\sigma,G}(H)$. To prove the converse, let $g\in
N_{\sigma,G}(H)$. Write $g=kp$ as in Lemma \ref{L:basic}. Then for
every $h_1\in H$, $h_2=\tau_g(h_1)\in H\subset K$. So
$\tau_k(h_1)=h_2\in H$, $\tau_p(h_1)=h_1$, that is, $k\in
N_{\sigma,K}(H), p\in Z_{\sigma,G}(H)$. Hence $g\in
N_{\sigma,K}(H)\cdot Z_{\sigma,G}(H)$. This proves
$N_{\sigma,G}(H)\subset N_{\sigma,K}(H)\cdot Z_{\sigma,G}(H)$.
\end{proof}

\begin{corollary}\label{C:normalizer}
Let $G$ be a connected Lie group with a maximal compact group $K$,
$H$ a closed subgroup of $K$. Then $N_G(H)=N_K(H)\cdot
Z_G(H)$.\qed
\end{corollary}

\begin{proposition}\label{P:same}
Let $G$ be a connected Lie group with an automorphism $\sigma$ of
finite order. Then\\
(i) For every maximal compact torus $T$ of $G_0^\sigma$, there
exists a $\sigma$-invariant maximal compact subgroup $K$ of $G$
such that $T$ is a maximal torus of $K_0^\sigma$;\\
(ii) For every $\sigma$-invariant maximal compact subgroup $K$ of
$G$ and every maximal torus $T$ of $K_0^\sigma$, $T$ is a maximal
compact torus of $G_0^\sigma$.
\end{proposition}

\begin{proof}
(i) Suppose $\sigma$ is of order $n$. Then $\Zn$ acts on $G$ by
$(m,g)\mapsto \sigma^m(g)$. Consider the semidirect product
$G\rtimes\Zn$. Then $T\rtimes\Zn$ is a compact subgroup of
$G\rtimes\Zn$. Let $\widetilde{K}$ be a maximal compact subgroup
of $G\rtimes\Zn$ containing $T\rtimes\Zn$, and let $K$ be the
identity component of $\widetilde{K}$. Then $K$ is a maximal
compact subgroup of $G$ containing $T$. Since $T\subset
G_0^\sigma$, $T\subset K_0^\sigma$. Moreover, since $T$ is a
maximal compact torus of $G_0^\sigma$, $T$ is a maximal torus of
$K_0^\sigma$.

(ii) Suppose $K$ is a $\sigma$-invariant maximal compact subgroup
of $G$, $T$ is a maximal torus of $K_0^\sigma$. It is obvious that
$T$ is a compact torus of $G_0^\sigma$. Let $T'$ be a maximal
compact torus of $G_0^\sigma$ containing $T$. We want to prove
that $T'=T$.

By (i), there is a $\sigma$-invariant maximal compact subgroup
$K'$ of $G$ such that $T'$ is a maximal torus of $(K')_0^\sigma$.
By \cite{Ho}, Chapter XV, Theorem 3.1, there exists $g\in G$ such
that $gK'g^{-1}=K$. Then $gT'g^{-1}$ is a maximal torus of
$g(K')_0^\sigma g^{-1}$. But for $k\in K$, $k\in g(K')^\sigma
g^{-1}\Leftrightarrow \sigma(g^{-1}kg)=g^{-1}kg\Leftrightarrow
k\in G^{\Inn(g\sigma(g)^{-1})\circ\sigma}$, where
$\Inn(g\sigma(g)^{-1})$ is the inner automorphism of $G$ induced
by $g\sigma(g)^{-1}$. Note that
$\Inn(g\sigma(g)^{-1})(K)=g\sigma(g)^{-1}K\sigma(g)g^{-1}=g\sigma(g^{-1}Kg)g^{-1}=gK'g^{-1}=K$,
that is, $K$ is $\Inn(g\sigma(g)^{-1})$-invariant. So
$g(K')^\sigma g^{-1}=K^{\Inn(g\sigma(g)^{-1})\circ\sigma}$, and
then $g(K')_0^\sigma
g^{-1}=K_0^{\Inn(g\sigma(g)^{-1})\circ\sigma}$. By Corollary
\ref{C:normalizer}, $g\sigma(g)^{-1}\in N_G(K)=K\cdot Z_G(K)$. So
$\Inn(g\sigma(g)^{-1})|_K$ is an inner automorphism of $K$. By
Theorem \ref{T:application}, $\rank K_0^\sigma=\rank
K_0^{\Inn(g\sigma(g)^{-1})\circ\sigma}$, that is, $\dim T=\dim
gT'g^{-1}=\dim T'$. Hence $T'=T$.
\end{proof}

Let $G$ be a connected Lie group with an automorphism $\sigma$ of
finite order, and let $T$ be a maximal compact torus of
$G_0^\sigma$. By Proposition \ref{P:same}, we can choose a
$\sigma$-invariant maximal compact subgroup $K$ of $G$ containing
$T$. The natural inclusion $N_{\sigma,K}(T)\hookrightarrow
N_{\sigma,G}(T)$ induces a natural map
$W(K,\sigma|_K,T)\rightarrow W(G,\sigma,T)$.

\begin{proposition}\label{P:iso}
Under the above conditions, the natural map
$W(K,\sigma|_K,T)\rightarrow W(G,\sigma,T)$ is an isomorphism.
\end{proposition}

\begin{proof}
By Lemma \ref{L:normalizer}, we have
\begin{align*}
&W(G,\sigma,T)\\
=&N_{\sigma,G}(T)/Z_{\sigma,G}(T)\\
=&N_{\sigma,K}(T)\cdot
Z_{\sigma,G}(T)/Z_{\sigma,G}(T)\\
\cong &N_{\sigma,K}(T)/N_{\sigma,K}(T)\cap
Z_{\sigma,G}(T)\\
=&N_{\sigma,K}(T)/Z_{\sigma,K}(T)\\
=&W(K,\sigma|_K,T).
\end{align*}
It is obvious that the inverse of the above isomorphism coincides
with the natural map $W(K,\sigma|_K,T)\rightarrow W(G,\sigma,T)$.
This proves the proposition.
\end{proof}

\begin{proposition}\label{P:finite}
$W(G,\sigma,T)$ is a finite group.
\end{proposition}

\begin{proof}
First assume $G$ is compact. Since $Z_\sigma(T)$ and $N_\sigma(T)$
are closed subgroups of $G$, they have finitely many connected
components. But by (ii) of Proposition \ref{P:lie algebra}, their
identity components coincide. So
$W(G,\sigma,T)=N_\sigma(T)/Z_\sigma(T)$ is finite.

Now suppose $G$ is noncompact. By Proposition \ref{P:same}, there
exists a $\sigma$-invariant maximal compact subgroup $K$ of $G$
such that $T$ is a maximal torus of $K_0^\sigma$. By Proposition
\ref{P:iso}, $W(G,\sigma,T)\cong W(K,\sigma|_K,T)$. But we have
proved that $W(K,\sigma|_K,T)$ is finite. So $W(G,\sigma,T)$ is
finite.
\end{proof}

The twisted conjugate action $\tau$ associated with $\sigma$
induces naturally an action of $W$ on $T$, defined by
$$w.t=\tau_g(t),$$ where $w\in W(G,\sigma,T), t\in T, g\in w$.
If $\sigma$ is of finite order and $n$ is a positive integer which
is divisible by the order of $\sigma$, then the finite subgroup
$E_n(T)=\{t\in T|t^n=e\}$ of $T$ is invariant under the natural
action of $W$ on $T$. In fact, if $t\in E_n(t)$, $g\in
N_\sigma(T)$, then
$\tau_g(t)^n=(gt\sigma(g)^{-1})^n=(gt\sigma(g^{-1}))
\sigma(gt\sigma(g^{-1}))\cdots\sigma^{n-1}(gt\sigma(g^{-1}))
=gt^n\sigma^n(g^{-1})=e$, that is, $\tau_g(t)\in E_n(T)$. So $W$
acts naturally on $E_n(T)$.

The following result is important for us to prove the main result
in this paper.

\begin{theorem}\label{T:conjugate}
Let $G$ be a connected compact Lie group with a $1$-semisimple
automorphism $\sigma$, $T$ a maximal torus of $G_0^\sigma$,
$W=W(G,\sigma,T)$. Then two elements $t_1,t_2\in T$ are
$\sigma$-conjugate if and only if they lie in the same $W$-orbit.
\end{theorem}

\begin{proof}
The ``if" part is obvious. To prove the converse, assume that
$t_1,t_2\in T$ are $\sigma$-conjugate, that is, there exists $g\in
G$ such that $t_2=\tau_g(t_1)$. Let $G_{t_1}=\{h\in
G|\tau_h(t_1)=t_1\}$, and let $(G_{t_1})_0$ be the identity
component of $G_{t_1}$. It is obvious that $T\subset (G_{t_1})_0$.
Note that $h\in G_{t_1}\Leftrightarrow
ht_1\sigma(h)^{-1}=t_1\Leftrightarrow
h=t_1\sigma(h)t_1^{-1}\Leftrightarrow h\in
G^{\Inn(t_1)\circ\sigma}$, that is,
$G_{t_1}=G^{\Inn(t_1)\circ\sigma}$. By Theorem
\ref{T:application}, $\rank G_0^\sigma=\rank (G_{t_1})_0$. So $T$
is a maximal torus of $(G_{t_1})_0$.

We claim that $g^{-1}Tg\subset(G_{t_1})_0$. In fact, for every
$t\in T$, we have
$\tau_{g^{-1}tg}(t_1)=\tau_{g^{-1}}\circ\tau_t\circ\tau_g(t_1)
=\tau_{g^{-1}}\circ\tau_t(t_2)=\tau_{g^{-1}}(t_2)=t_1$. So
$g^{-1}Tg\subset G_{t_1}$. But $g^{-1}Tg$ is connected, so
$g^{-1}Tg\subset(G_{t_1})_0$.

Since $T$ and $g^{-1}Tg$ are maximal tori of $(G_{t_1})_0$, there
exists $x\in (G_{t_1})_0$ such that $xTx^{-1}=g^{-1}Tg$. Let
$y=gx$, then $y^{-1}Ty=T$, and then
$\tau_{y^{-1}}(T)=y^{-1}T\sigma(y)=y^{-1}Tyy^{-1}t_2\sigma(y)
=Tx^{-1}g^{-1}t_2\sigma(g)\sigma(x)=T\tau_{x^{-1}}(t_1)=Tt_1=T$.
So $y\in N_\sigma(T)$. But
$\tau_y(t_1)=\tau_g\circ\tau_x(t_1)=\tau_g(t_1)=t_2$. So $t_1$ and
$t_2$ lie in the same $W$-orbit. This completes the proof of the
theorem.
\end{proof}

%%%%%%%%%%%%%%%%%%%%%%%%%%%%%%%%%%%%%%%%%%%%%%%%%%%%%%%%%%% Section 4

\section{Twisted Weyl groups and nonabelian cohomology}

In this section we define the twisted Weyl group associated with
an $A$-module structure on a Lie group, and give the proof Theorem
\ref{T:main}.

Let $A$ be a cyclic group, that is, $A\cong\ZZ$ or $A\cong\Zn$,
and let $G$ be a connected Lie group with an $A$-module structure.
\emph{If $A\cong\ZZ$, we always assume that $G$ is compact and the
$\ZZ$-module structure on $G$ is $1$-semisimple}, that is, the
action of a generator of $\ZZ$ on $G$ is $1$-semisimple. Let $T$
be a maximal compact torus of $G^A$. Then for every generator
$\sigma$ of $A$, we can construct the twisted Weyl group
$W(G,\sigma,T)$ associated with $\sigma$.

\begin{proposition}\label{P:module}
Under the above conditions, $W(G,\sigma,T)$ is independent of the
choice of the generator $\sigma$ of $A$.
\end{proposition}

\begin{proof}
It is sufficient to prove that $Z_\sigma(T)$ and $N_\sigma(T)$ are
independent of the choice of $\sigma$. By (i) of Proposition
\ref{P:lie algebra}, $Z_\sigma(T)\subset G^A$. So
$Z_\sigma(T)=Z_{G^A}(T)$, which is obviously independent of the
choice of $\sigma$.

To prove that $N_\sigma(T)$ is independent of the choice of
$\sigma$, let $\sigma'$ be a generator of $A$ different from
$\sigma$. If $A\cong\ZZ$, then $\sigma'=\sigma^{-1}$. Let $g\in
N_\sigma(T)$, then
\begin{align*}
&gT\sigma'(g)^{-1}=gT\sigma^{-1}(g^{-1})\\
=&\sigma^{-1}(\sigma(g)Tg^{-1})=\sigma^{-1}((gT\sigma(g)^{-1})^{-1}))\\
=&\sigma^{-1}(T^{-1})=T.
\end{align*}
So $g\in N_{\sigma'}(T)$. Hence $N_\sigma(T)\subset
N_{\sigma'}(T)$. By symmetry, we have
$N_\sigma(T)=N_{\sigma'}(T)$. If $A\cong\Zn$, then there exist
positive integers $r$ and $s$ such that $\sigma'=\sigma^r$,
$\sigma=(\sigma')^s$. Let $g\in N_\sigma(T)$, then
\begin{align*}
&gT\sigma'(g)^{-1}=gT\sigma^r(g^{-1})\\
=&(gT\sigma(g^{-1}))(\sigma(g)T\sigma^2(g^{-1}))\cdots(\sigma^{r-1}(g)T\sigma^r(g^{-1}))\\
=&(gT\sigma(g^{-1}))\sigma(gT\sigma(g^{-1}))\cdots\sigma^{r-1}(gT\sigma(g^{-1}))\\
=&T^r=T.
\end{align*}
Hence $N_\sigma(T)\subset N_{\sigma'}(T)$. By symmetry, we also
have $N_\sigma(T)=N_{\sigma'}(T)$.
\end{proof}

In virtue of Proposition \ref{P:module}, we can define the
\emph{twisted Weyl group associated with the $A$-module $G$} by
$$W(G,A,T)=W(G,\sigma,T).$$ By Proposition \ref{P:independent},
$W(G,A,T)$, as an abstract group, is independent of the choice of
$T$. We will simply denote $W(G,A,T)$ by $W$ if the omitted data
are explicit from the context.

Now consider the cohomology $H^1(A,T)$. Since $T$ is abelian and
$A$ acts trivially on $T$, $H^1(A,T)$ coincides with the set of
cocycles $Z^1(A,T)$, and an element $\alpha\in Z^1(A,T)$ is a
homomorphism $\alpha:A\rightarrow T$. Then $W(G,A,T)$ acts
naturally on $Z^1(A,T)\cong H^1(A,T)$ by
$$(w.\alpha)(a)=g\alpha(a)a(g)^{-1},$$ where $w\in W(G,A,T),
\alpha\in Z^1(A,T), a\in A, g\in w$.

Since $A$ is cyclic, a homomorphism $\alpha:A\rightarrow T$ is
determined by its value $\alpha(\sigma)$ on a generator $\sigma$
of $A$. So if a generator $\sigma$ of $A$ is chosen, we have
natural identifications $T\cong Z^1(\ZZ,T)\cong H^1(\ZZ,T)$ and
$E_n(T)\cong Z^1(\Zn,T)\cong H^1(\Zn,T)$. Under these
identifications, it is obvious that the natural action of
$W(G,\ZZ,T)$ on $H^1(\ZZ,T)$ coincides with the natural action of
$W(G,\sigma,T)$ on $T$, and the natural action of $W(G,\Zn,T)$ on
$H^1(\Zn,T)$ coincides with the natural action of $W(G,\sigma,T)$
on $E_n(T)$.

The natural inclusion $T\hookrightarrow G$ induces a natural map
$H^1(\ZZ,T)\rightarrow H^1(\ZZ,G)$, which obviously reduces to a
map $W\backslash H^1(\ZZ,T)\rightarrow H^1(\ZZ,G)$. Now we are
prepared to prove the main result of this paper.

\begin{theorem}\label{T:main'}
Let $A$ be a cyclic group, $G$ a connected Lie group with an
$A$-module structure. If $A\cong\ZZ$, we assume that $G$ is
compact and the $\ZZ$-module structure on $G$ is $1$-semisimple.
Let $T$ be a maximal compact torus of $G_0^A$, $W=W(G,A,T)$ the
associated twisted Weyl group. Then the map $W\backslash
H^1(A,T)\rightarrow H^1(A,G)$ is a bijection.
\end{theorem}

\begin{proof}
By Theorems 4.1 and 5.1 in \cite{AW}, the natural map
$H^1(A,T)\rightarrow H^1(A,G)$ is surjective, so $W\backslash
H^1(A,T)\rightarrow H^1(A,G)$ is surjective.

To prove the injectivity, we first assume $A\cong\ZZ$. Choose a
generator $\sigma$ of $\ZZ$. Then there is an identification
$H^1(\ZZ,T)\cong T$, under which the natural action of
$W(G,\ZZ,T)$ on $H^1(\ZZ,T)$ coincides with the natural action of
$W(G,\sigma,T)$ on $T$. Suppose $t_1, t_2\in T\cong H^1(\ZZ,T)$
have the same image in $H^1(\ZZ,G)$ under the natural map
$H^1(\ZZ,T)\rightarrow H^1(\ZZ,G)$. Then $t_1$ and $t_2$ are
$\sigma$-conjugate. By Theorem \ref{T:conjugate}, $t_1$ and $t_2$
lie in the same $W(G,\sigma,T)$-orbit. So $W\backslash
H^1(\ZZ,T)\rightarrow H^1(\ZZ,G)$ is injective.

Now assume that $A\cong\Zn$. By Proposition \ref{P:same}, there
exists a maximal compact subgroup $K$ of $G$ which is a
$\Zn$-submodule such that $T$ is a maximal torus of $K_0^{\Zn}$.
By Theorem 3.1 in \cite{AW}, the natural map
$H^1(\Zn,K)\rightarrow H^1(\Zn,G)$ is bijective. So it is
sufficient to prove that $W\backslash H^1(\Zn,T)\rightarrow
H^1(\Zn,K)$ is injective. Let $\sigma'$ be a generator of $\Zn$.
Then there is an identification $H^1(\Zn,T)\cong E_n(T)$. Suppose
$t_1, t_2\in E_n(T)\cong H^1(\Zn,T)$ have the same image in
$H^1(\Zn,K)$ under the natural map $H^1(\Zn,T)\rightarrow
H^1(\Zn,K)$. Then $t_1$ and $t_2$ are $\sigma'|_K$-conjugate. By
Theorem \ref{T:conjugate}, $t_1$ and $t_2$ lie in the same
$W(K,\sigma'|_K,T)$-orbit. But by Proposition \ref{P:iso},
$W(G,\sigma',T)\cong W(K,\sigma'|_K,T)$, and it is obvious that
the natural actions of $W(G,\sigma',T)$ and $(K,\sigma'|_K,T)$ on
$T$ coincide. So $W\backslash H^1(\ZZ,T)\rightarrow H^1(\ZZ,G)$ is
injective.
\end{proof}

\begin{remark}
By Proposition \ref{P:same}, for a connected Lie group $G$ with a
$\Zn$-module structure, the family $\{T|T \ \text{is a maximal
compact torus of} \ G_0^{\Zn}\}$ coincides with the family
$\bigcup_K\{T|T \ \text{is a maximal torus of} \ K_0^{\Zn}\}$,
where $K$ runs through all maximal compact subgroups of $G$ which
are $\Zn$-submodules.
\end{remark}

\end{document}